
\documentstyle{amsppt}
\magnification=\magstep 1
\TagsOnRight
\NoBlackBoxes
\leftheadtext{V\.G\. Pestov}
\rightheadtext{Banach analytic spaces}
\def\norm #1{{\left\Vert\,#1\,\right\Vert}}
\def\BR #1{{B_{\R}(#1)}}
\def\BC #1{{B_{\C}(#1)}}
\def\B #1{{B(#1)}}
\def\R {{\Bbb R}}
\def\C {{\Bbb C}}

\def\e{{\epsilon}}

\def\supp{{\operatorname{supp}\,}}
\def\span{{\operatorname{span}\,}}
\def\Card{{\operatorname{Card}\,}}
\def\Lip{{\operatorname{Lip}\,}}
\def\dist{{\operatorname{dist}\,}}
\def\diam{{\operatorname{diam}\,}}
\def\QED{\nobreak\quad\ifmmode\roman{Q.E.D.}\else{\rm Q.E.D.}\fi}
\def\ml #1\endml{\email #1\endemail}
\long\def\block #1\endblock{\vskip 6pt
	{\leftskip=3pc \rightskip=\leftskip
	\noindent #1\endgraf}\vskip 6pt}
\long\def\ext #1\endext{\block #1\endblock}

\def\ml #1\endml{\email #1\endemail}
\topmatter
\title
Douady's conjecture on Banach analytic spaces
$\dag$
\endtitle
\author
Vladimir G. Pestov
\endauthor
\address Department of Mathematics, Victoria University of Wellington,
P.O. Box 600, Wellington, New Zealand
\endaddress
\ml vladimir.pestov\@vuw.ac.nz\endml
\date June 20, 1994\enddate

\abstract{We show that, as conjectured by Adrien Douady back in 1972,
every complete metric space is homeomorphic (moreover, isometric)
to the locus of zeros of an analytic map between two
Banach spaces.}
\endabstract
\subjclass{Primary 32K05; Secondary 16W30, 32C18, 46B04}
\endsubjclass
\keywords{Banach analytic spaces, free Banach spaces, Banach coalgebras}
\endkeywords
\endtopmatter
\document
\footnote""{$\dag$ Research Report RP-94-144, Department of Mathematics,
Victoria University of Wellington, June 1994.}

\subheading{1. Introduction}
The concept of a complex analytic space \cite{Fi}
admits infinite-dimensional generalizations,
but their geometric content is not so clear.
Banach analytic spaces as defined by Douady
\cite{D1} are topological spaces endowed
with a sheaf of $B$-valued maps for every Banach space $B$
in a functorial way and locally isomorphic with models
(loci of zeros of analytic maps between complex Banach spaces).
In his 1972 note \cite{D2}, Douady had shown that every compact
metric space is homeomorphic to the locus of zeros of an analytic map
between two Banach spaces, and therefore supports the structure of a
Banach analytic space.
He has conjectured that it must be true for every
complete metric space, and argued that
``there is no point in studying Banach analytic spaces for their
own sake.''

\par Here we prove Douady's conjecture in a somewhat stronger form.

\proclaim{Main Theorem} A complete metric space is isometric
to the locus of zeros of an analytic map between two complex
Banach spaces.
\endproclaim

\proclaim{Corollary} A paracompact topological space admits the
structure of a Banach analytic space if and only if it is metrizable
with a complete metric.
\endproclaim

We sharpen the original argument for compact spaces
by dualizing it and considering instead of the algebra of
Lipschitz functions on a metric space $X$ its
predual, the so-called free Banach space.
In the process, multiplication is replaced
with comultiplication, and a
complete metric space of diameter $\leq 2$
embeds as a subspace of group-like
elements in a countial Banach coalgebra. From this result, we
are able to deduce the Main Theorem.

\subheading{2. Free Banach spaces}
 The construction of a free normed space,
going back to Arens and Eells \cite{AE},
is relatively well understood \cite{M, Fl1, Fl2, R, P}.
Completing a free normed space, one obtains a free
Banach space \cite{P}.

\definition{Definition 1} \cite{Fl1,Fl2}
A {\it norm pair} $(\rho,\alpha)$ on a set $X$ consists of a
metric $\rho$ and a non-negative
real-valued function $\alpha$ such that
$$\vert\alpha(x)-\alpha(y)\vert\leq\rho(x,y)\leq
\alpha(x)+\alpha(y)$$
for all $x,y\in X$.
\enddefinition

In absence of a better name, we refer to a triple
$(X,\rho,\alpha)$, where $(\rho,\alpha)$ is a norm pair on $X$,
as a {\it normed set}. A mapping $f\colon X\to Y$
between normed sets is {\it Lipschitz with constant} $L$ if
for every $x,y,z\in X$ one has $\rho_Y(f(x),f(y))\leq
L\rho_X(x,y)$ and $\alpha_Y(z)\leq L\alpha_X(z)$.
A Lipschitz map with constant $1$ is a {\it contraction} of
normed spaces. (Though the category of normed spaces and contractions
is naturally equivalent to the category of pointed metric spaces and
base-preserving contractions, the former one is more convenient for
our purposes.)
A subset of a normed set naturally inherits a normed subset
structure.
Every normed space $E$ becomes a normed set if one sets
$\rho(x,y)=\norm{x-y}$ and $\alpha(z)=\norm z$.
The following is a ``completed'' version of Flood's theorem
\cite{Fl1,Fl2}.

\proclaim{Theorem 1}
Let $X=(X,\rho_X,\alpha_X)$ be a normed set. There exists a complete
normed space $B(X)$ and an embedding of $X$ into
$\B X$ as a normed subset
such that every contraction $f$ from $X$ to a complete normed space
$E$ lifts to a unique linear contraction $\bar f\colon B(X)\to E$.
The pair consisting of $B(X)$ and embedding $X\hookrightarrow B(X)$
is essentially unique. Elements of $X$ are linearly independent.
\qed\endproclaim

The following is proved by an obvious renormalization of $E$.

\proclaim{Assertion 1} Let $X$ be a normed set and
let $f$ be a Lipschitz map with constant
$L$ from $X$ to a complete normed space
$E$. Then $f$ extends to a unique
linear map $\bar f\colon B(X)\to E$ with $\norm{\bar f}\leq L$.
\qed\endproclaim

Sometimes we have to distinguish between the real and complex
versions of the free Banach space, $\BR X$ and $\BC X$ respectively.

\proclaim{Assertion 2} {\rm (\cite{Fl1, Fl2}; cf. \cite{R}.)
Let $x$ be in the linear span of a normed set
$X$ in $\BR X$. Then for some finite
collection $\lambda_i,\mu_j\in\R$ and
$x_i,y_i,z_j\in\supp x$, one has
$$x=\sum_i\lambda_i(x_i-y_i)+\sum_j\mu_jz_j$$
and
$$\norm x = \sum_i\vert\lambda_i\vert\rho_X(x_i-y_i)+
\sum_j\vert\mu_j\vert\alpha_X(z_j)$$
\qed\endproclaim

\proclaim{Assertion 3} Let $X$ be a normed subset of $Y$.
Then the embedding $X\to Y$ extends to an
isometric embedding of Banach spaces $\BR X\to \BR Y$.
\endproclaim

\demo{Proof} Assertion 2 implies
that the canonical linear contraction $\BR X\to \BR Y$
preserves the norm of every
element from $\span X$. The rest is clear.
\qed\enddemo

A normed set $Y$ is {\it injective} if for every normed set
$M$ and a normed subset $N\subseteq M$,
an arbitrary contraction $f\colon N\to Y$ extends to a
contraction $\tilde f\colon M\to Y$.

\proclaim{Assertion 4} {\rm (Cf. \cite{N}.)}
Every normed set $X$ embeds in an injective normed set $Y$ such
that the metric $\rho_Y$ is complete.
\endproclaim

\demo{Proof} As any normed space, $\BR X$ embeds into a
complete normed space $E$ with the $1$-extension property
(\cite{L}, \S 11, the proof of Coroll. 1).
The space $E$ viewed as a normed set has the desired
properties of $Y$.
It is enough to apply the $1$-extension property
of $E$ to the pair $(\BR N,\BR M)$, where $(N,M)$ are as
above, keeping in mind the universal property of the
free Banach space and Assertion 3.
\qed\enddemo

\proclaim{Theorem 2} A normed set $X$
is weakly closed in $\B X$ if and only if the metric $\rho_X$ is
complete.
\endproclaim

\demo{Proof} Necessity being obvious, we
concentrate on sufficiency. Applying Assertions 4 and 3,
one can assume $X$ to be injective.
Let $x\in \B X\setminus X$. Since $X$ is complete,
one must have $\dist(x,X)>0$.
Fix an $y\in\span X$ with
$\norm{x-y}<\dist(x,X)/2$; clearly, $\dist(y,X)>\dist(x,X)/2$.
For every $a\in X$ and $z\in\{0\}\cup\supp y$, set
$$f(a)_z=\min\{\norm{z-a}-\norm z,\diam(\{0\}\cup\supp y)\}$$
and
$$f(a)=(f(a)_z)_{z\in\{0\}\cup\supp y}\in l_{\infty}^k$$
where $k=\Card(\supp y)+1$.

The mapping $f\colon X\to K=\overline{f(X)}$ is a contraction
relative the structure of a normed set induced on a (compact) subset
$K$ from $l_{\infty}^k$.
Extend $f$ to a linear contraction
$\bar f\colon\B X\to\B{K}$. We claim that
$\dist(\bar f(y),K)=\dist(y,X)$.
While the inequality $\leq$ is obvious, assume that
there exists a $b\in K$
with $\norm{\bar f(y)-b}<\dist(y,X)$.
Since the restriction of $f$ to $\{0\}\cup\supp y$ is an
isometric embedding of normed pairs
and $X$ is injective, there is a contraction
$\phi\colon K\to X$ with $\phi\circ f\vert_{\{0\}\cup\supp y}=
\operatorname{Id}_{\{0\}\cup\supp y}$. Now for the extension
$\bar\phi\colon \B{K}\to\B X$ one has
$\norm{y-\phi(b)}=\norm{\bar\phi(\bar f(y)-b)}
\leq \norm{\bar f(y)-b}<\dist(y,X)$,
a contradiction since $\phi(b)\in X$.

As a consequence of the established equality,
$\dist(\bar f(x),K)\geq\dist(\bar f(y),K)-\norm{x-y}>0$ and
$\bar f(x)\notin K$. Finally,
$K$ is compact and therefore weakly
closed in $\B K$, and $\bar f(X)\subseteq K$.
\qed\enddemo

\proclaim{Lemma 1} Let $X$ be a normed set such that the metric
$\rho_X$ is complete and $\alpha_X\geq 1$. Let $x\in\B X\setminus X$.
Then there exists a contraction $f$ from $X$ to a normed
set $K$ with $\alpha_K\equiv 1$ such that $\bar f(x)\in\B K\setminus K$.
\endproclaim

\demo{Proof}  There is a
linear contraction $f$ from $\B X$ to a finite-dimensional
normed space $E$ with $x\notin\overline{f(X)}$ (Theorem 2).
Make a compact subspace $K=\overline{f(X)}$ of $E$ into
a normed set by letting $\rho_K(x,y)=L^{-1}\norm{x-y}$ and
$\alpha_K\equiv 1$, where $L=\max\{1,\diam K/2\}$. The
map $f\colon X\to K$ is easily verified to be a contraction of
normed sets and
therefore extends to a linear contraction
$\bar f\colon\B X\to\B K$.
As an application of Assertion 1, $f$ factors through $\B K$, thence
$\bar f(x)\notin K$.
\qed\enddemo

\proclaim{Lemma 2} Let $(X,\rho)$ be a complete metric space
and let $x_0\in X$.
Set $\alpha=1+\rho(x,x_0)$, $d(x,y)=\min\{\rho(x,y),1\}$, and
$\beta\equiv 1$. Then the linear contraction
$i\colon \B {X,\rho,\alpha}\to \B {X,d,\beta}$
extending the identity mapping
$\operatorname{Id}_X$ has the property
$i^{-1}(X)=X$.
\endproclaim

\demo{Proof} Let $x\in\B {X,\rho,\alpha}\setminus X$.
Choose $f$ and $K$ as in Lemma 1.
The map $f$ is easily checked to be a Lipschitz map with constant
$2$ from
$(X,d,\beta)$ to $K$ and therefore it extends to a bounded linear map
$\tilde f\colon \B {X,d,\beta}\to\B K$. Clearly,
$\bar f=\tilde f\circ i$
and therefore $i(x)\notin X$, as desired.
\qed\enddemo

\subheading{3. Banach coalgebras}
\definition{Definition 2} We call a pair
$C=(C,\Delta)$ consisting of a Banach space $C$ and
a bounded linear map $\Delta\colon C\to C\hat\otimes C$
({\it comultiplication})
a {\it Banach coalgebra}. A Banach coalgebra $C$ is
{\it counital} if it has a {\it counit}, that is, a
bounded linear mapping $\e\colon C\to \C$ with
$$(\operatorname{Id}_C\hat\otimes\epsilon)\circ\Delta
=(\epsilon\hat\otimes\operatorname{Id}_C)\circ\Delta =
\operatorname{Id}_C$$
(Here $\hat\otimes$ stands for the projective tensor
product.)
\enddefinition

\par
We call an element $x$ of a Banach coalgebra $C$
{\it group-like} (cf. \cite{A}, 2.1.1) if
$$\Delta x=x\otimes x~~\text{and}~~\e(x)=1$$

\proclaim{Assertion 5}
The subspace of group-like elements of a counital
Banach coalgebra $C$
forms the locus of zeros of an analytic map from $C$ to a
Banach space.
\endproclaim

\demo{Proof} The map $x\mapsto x\otimes x$ from $C$ to
$C\hat\otimes C$ is a homogenuous
polynomial map of degree $2$, which is bounded on the unit ball
and therefore continuous.
(Cf. Prop. 1 in \cite{D1}.)
Therefore, the map
$$C\ni x\mapsto (x\otimes x -\Delta x, \e(x)-1_{\C})\in
(C\hat\otimes C)\oplus\C$$
is analytic (moreover, continuous polynomial of degree $2$).
Zeros of this map are exactly the group-likes of $C$.
\qed\enddemo

\proclaim{Theorem 3} Suppose a metric space $X$ of
diameter $\leq 2$ is complete. Make $X$ into a normed set
by letting $\alpha_X\equiv 1$.
Then $\BC X$ admits a unique structure of a counital Banach coalgebra
such that the subspace of group-like elements coincides
with $X$.
\endproclaim

\demo{Proof}
Assertion 1 implies the existence of a comultiplication
(the bounded linear operator extending
to $\BC X$ the Lipschitz map $X\ni x\mapsto x\otimes x\in
\BC X\hat\otimes\BC X$)
and a counit (the extention of the contraction
$X\ni x\mapsto 1\in\C$). Uniqueness of these
coalgebra operations is clear.

Let $x\in\BC X\setminus X$. Choose a contraction $f$ as in
Lemma 1 and notice that $\BC K$ also becomes a Banach coalgebra
as above.
Definitely, $\bar f$ preserves
the comultiplication and commutes with the counit map.
The Banach space dual to $\BC K$ is the space of all complex-valued
Lipschitz functions, $\Lip(K)$, on $K$ with the usual norm topology.
Denote by $\langle\,,\,\rangle$ the
natural pairing between $\BC K\hat\otimes\BC K$ and the space of
all bounded bilinear functionals on $\BC K\times\BC K$.
Every pair $g,h\in \Lip(K)$ determines such a functional
$g\otimes h\colon(a,b)\mapsto \bar f(a)\bar g(b)$.
It is easy to see that
$\langle \Delta x,g\otimes h \rangle = \overline{(fg)}(x)$ and
$\langle x\otimes x,g\otimes h \rangle = \bar f (x)\bar g(x)$.
Assuming that $x$ is a group-like in $\B X$,
its image $\bar f(x)$ is a group-like in
$\B K$; therefore $\overline{(fg)}(x)=
\bar f (x)\bar g(x)$ and $x$ is a multiplicative
functional on the Banach algebra $\Lip(K)$.
(Remark that also $x(1)=1$.)
However, the only multiplicative functionals on the algebra
of Lipschitz functions on a compact metric space $K$
are known to be (defined by) the elements of $K$.
(Cf. \cite{D2};\cite{G}, ex. I.3.4.)
\qed\enddemo

\subheading{4. Proofs} {\it Main Theorem:} If the diameter $L>0$
of $X$ does not
exceed $2$, Theorem 3 and Assertion 5 do the job. Otherwise,
$X$ is the locus of zeros of the composition of a linear
contraction $i\colon \B {X,\rho,\alpha}\to \B {X,d,\beta}$
from Lemma 2 and the polynomial map from Assertion 5. \qed
\smallskip
\noindent{\it Corollary:}
It is enough to show that a paracompact space $X$
locally homeomorphic to a complete metric space is itself
metrizable with a complete metric. Choose a locally finite open
cover of $X$, $\gamma=\{U_{\alpha}\colon\alpha\in A\}$, and complete
metrics $\rho_{\alpha}\leq 1$ on $\overline{U_{\alpha}}$.
One can assume that each $U_\alpha$ is a unit ball of $\rho_\alpha$.
Endow the quotient space
$X_\alpha=X/(X\setminus U_\alpha)$
with the largest
metric making the quotient map $\pi_\alpha\colon X\to X_\alpha$
a contraction. The metric $\rho(x,y)=\sup_{\alpha}\rho_\alpha
(\pi_\alpha x,\pi_\alpha y)$ on
$\prod_{\alpha}X_\alpha$'s is complete and the map
$\pi(x)=(\pi_\alpha(x_\alpha))_\alpha$
is a homeomorphic embedding with a closed range. \qed

\subheading{5. Open question}
Every topological space admitting the structure of a Banach
analytic space is locally homeomorphic to a complete
metric space. Is the inverse correct?

\enddocument
\bye